\documentclass[aop,MSNbibl,secthm,dvips]{arximspdf}

\doi{10.1214/09-AOP491}
\volume{38}
\issue{2}
\pubyear{2010}
\firstpage{632}
\lastpage{660}

\makeatletter
\DeclareMathAlphabet\mathcaligr{OMS}{cmsy}{m}{n}

\newtheorem{lem}{Lemma}[section]
\newtheorem{pro}{Proposition}[section]

\newproclaim{defn}{Definition}[section]
\newproclaim{qu}{Question}
\newproclaim{rema}{Remark}
\def\bsuffix #1{#1}
\makeatother

\begin{document}
\begin{frontmatter}

\title{The asymptotic shape theorem for generalized first passage percolation}
\runtitle{A general asymptotic shape theorem}

\begin{aug}
\author[A]{\fnms{Michael} \snm{Bj\"orklund}\ead[label=e1]{mickebj@math.kth.se}\corref{}}
\runauthor{M. Bj\"orklund}
\affiliation{Royal Institute of Technology}
\address[A]{Royal Institute of Technology (KTH)\\
Department of Mathematics\\
S-100 44 Stockholm\\
Sweden\\
\printead{e1}} 
\end{aug}

\received{\smonth{2} \syear{2009}}
\revised{\smonth{5} \syear{2009}}

%
\begin{abstract}
We generalize the asymptotic shape theorem in first passage percolation
on $\mathbb{Z}^d$ to cover the case of general semimetrics. We prove
a structure theorem for equivariant semimetrics on topological groups
and an extended version of the maximal inequality for
$\mathbb{Z}^d$-cocycles of Boivin and Derriennic in the vector-valued case.
This inequality will imply a very general form of Kingman's subadditive
ergodic theorem. For certain classes of generalized first passage
percolation, we prove further structure theorems and provide rates of
convergence
for the asymptotic shape theorem. We also establish a general form of
the multiplicative ergodic theorem of Karlsson and Ledrappier
for cocycles with values in separable Banach spaces with the
Radon--Nikodym property.
\end{abstract}

%

\end{frontmatter}

\section{Introduction}

First passage percolation was introduced by Hammersley and Welsh in
\cite{HaWe}. A detailed description of the model
is given in Section \ref{class}. The theory can be roughly described as
the study of the generic large-scale geometry of semimetric spaces,
where the semimetric is allowed to vary measurably. The classical case
deals with the space $\mathbb{Z}^d$ and semimetrics induced
by random weights on the edges of the standard Cayley graph of $\mathbb{Z}^d$.
However, the setup easily extends to general groups.

In this paper, we introduce the notion of a random semimetric. Let $G$
be a locally compact group and suppose that $G$ acts on a
probability space $(X,\mu)$, where $\mu$ is invariant under the action
of $G$. We say that the action is \textit{ergodic} if the invariant sets
are either
null or conull, and \textit{quasi-invariant} if it preserves the measure
class of $\mu$. Suppose that $(Y,\nu)$ is a $\sigma$-finite measure
space. A \textit{random semimetric} on~$Y$, modeled on the $G$-space
$X$, is
a map $\rho\dvtx X \times Y \times Y \rightarrow[0,\infty)$ such
that $\rho_x$
is a~semimetric for almost every $x$ in $X$ and
\[
\rho_{g.x}(y,y') = \rho_x(g.y,g.y')
\]
for all $y,y'$ in $Y$, $g$ in $G$ and $x \in X$, and for all $y,y'$ in
$Y$, the map
\[
x \mapsto\rho_x(y,y')
\]
is measurable. In general, these objects are very complicated and form
the basis of subadditive ergodic theory. However, it turns out that
all random semimetrics can be realized as norms of additive cocycles
with values in large Banach spaces. The definition of a Gelfand
cocycle is given in Section \ref{coho} and is rather technical, but
turns out to be useful, in view of the following theorem.
\begin{thm}[(Structure theorem)] \label{main1}
Let $G$ be a locally compact, second countable group. Suppose that
$(X,\mu)$ is a probability measure space
with a $G$-invariant ergodic measure $\mu$. Suppose that $(Y,\nu)$ is a
$G$-space with a quasi-invariant $\sigma$-finite
measure $\nu$. If $\rho$ is a random $G$-equivariant semimetric on $Y$,
modeled on the $G$-space $(X,\mu)$, then there exists a
Gelfand $L^{1}(Y,\nu)$-cocycle, with respect to the left-regular
representation of $G$ on $L^{\infty}(Y,\nu)$ on the $G$-space
$(X,\mu)$, such that
\[
\rho_x(y,y') = \Vert s_x(y,y')\Vert _{L^{\infty}(Y,\nu)}.
\]
\end{thm}

We will refer to a random semimetric $\rho$ on a space $Y$ as \textit
{generalized first passage percolation} on $Y$. In view of Theorem \ref{main1},
the study of generalized first passage percolation is equivalent to the
study of Gelfand cocycles with values in $L^{\infty}(Y,\nu)$. However,
any Gelfand
cocycle with values in the dual of a Banach space $B$ defines a random
semimetric. In Sections \ref{horo} and \ref{conv}, we will restrict
the class
of Banach spaces under consideration and this will allow us to
establish certain structure theorems which are not known for classical
first passage
percolation. For instance, we determine the horofunctions of random
semimetric spaces when the cocycles take values in separable Hilbert
spaces and
we prove an analog of Kesten's celebrated inequality for classical
first passage percolation in this context.

However, the main result of this paper is the following extension of
Boivin's asymptotic shape theorem to general random semimetrics.
\begin{thm}[(Asymptotic shape theorem)] \label{main2}
Suppose that $\rho$ is a random $\mathbb{Z}^d$-semimetric modeled on
an ergodic
$\mathbb{Z}^d$-space $(X,\mu)$, where $\mu$ is a probability measure.
Suppose that $\rho(0,n)$ is in $L^{d,1}(X,\mu)$ for every $n \in
\mathbb{Z}^d$.
There then exists a~seminorm $L$ on $\mathbb{R}^d$ such that
\[
\lim_{|n| \rightarrow\infty} \frac{\rho_x(0,n) - L(n)}{|n|} = 0
\]
almost everywhere on $(X,\mu)$.
\end{thm}

 This result was only known for a certain class of \textit{inner} random
semimetrics on $\mathbb{Z}^d$~\cite{Bo90}. It can be proven \cite
{De} that the
integrability
condition on the cocycle to belong to the Lorentz space $L^{d,1}(X)$ is
sharp. The unit ball of the semimetric $L$ roughly describes the
generic asymptotic
shape of large balls in $\mathbb{Z}^d$ with the random semimetric
$\rho$. In
the general ergodic situation, essentially all convex shapes can be
attained as
asymptotic shapes. This is a result of H\"aggstr\"om Meester
\cite{HaMe}.

\section{Generalized first passage percolation}

\subsection{Bochner--Lorentz spaces}

In the following sections, we will make use of certain classes of
function spaces introduced by Lorentz in \cite{Lo}. It is
straightforward to extend the definition to cover
the case of vector-valued functions, and we will do so. Before we give
the definition of the necessary function spaces, we recall some basic
notions and useful facts about
measurability of vector-valued functions. Let $B$ be a Banach space and
$(X,\mathfrak{F},\mu)$ a measure space. A \textit{simple} function
$f \dvtx X
\rightarrow B$ is a function
on the form
\[
f = \sum_{k=1}^{n} c_k \chi_{A_k},
\]
where $A_k$ are elements of $\mathfrak{F}$ and $c_k$ are elements in
$B$. A function $f \dvtx X \rightarrow B$ is \textit{Bochner measurable} (or
strongly measurable) if there exists a sequence of
simple functions $f_n \dvtx X \rightarrow B$ such that $\Vert f_n-f\Vert _B
\rightarrow0$. A
function $f \dvtx X \rightarrow B$ is \textit{weakly measurable} if
\[
x \mapsto\langle\lambda, f(x) \rangle
\]
is measurable for every $\lambda$ in $B^*$, where $B^*$ is the dual of
$B$. A function $f \dvtx X \rightarrow B^*$ is \textit{weak*-measurable} if
\[
x \mapsto\langle\lambda, f(x) \rangle
\]
is measurable for every $\lambda$ in $B$, canonically identified with
an element of $B^{**}$.

We now turn to the definition of the function spaces. Let $1 \leq p,q
\leq\infty$, and suppose that $f$ is a complex-valued measurable
function on $X$. We define
\[
f^{*}(t) = \inf\{ s > 0 | d_f(s) \leq t \},
\]
where $d_f$ is the distribution function of $f$, that is,
\[
d_f(\alpha) = \mu\bigl(\{ x \in X | |f(x)| > \alpha \}\bigr),\qquad \alpha
\geq0.
\]
We define the norm
\[
\Vert f\Vert _{p,q} =
\cases{
\biggl(\displaystyle\int_{0}^{\infty} (t^{1/p} f^{*}(t))^{q} \frac {dt}{t}\biggr)^{1/q},
&\quad
$\mbox{if $q < \infty$},$ \cr
\displaystyle\sup_{t > 0} t^{1/p} f^{*}(t), &\quad $\mbox{if $q=\infty$.}$}
\]
We denote the set of all $f$ with $\Vert f\Vert _{p,q} < \infty$ by
$L^{p,q}(X)$. With the above norm, this is a Banach space, usually
referred to as the \textit{Lorentz space} with indices $p$ and~$q$. For
instance, we see that $L^{p,p}(X) = L^{p}(X)$.

The extension to vector-valued functions is straightforward: we say
that a~weak*-measurable function $f \dvtx X \rightarrow B $ is in
$L^{p,q}_{w^*}(X,B^*)$ if there exists a nonnegative function $g$ on
$X$ with
finite $L^{p,q}(X)$-norm such that $\Vert f(x)\Vert _B \leq g(x)$ almost
everywhere. Note that $\Vert f\Vert _B$ is not necessarily measurable on $X$.
If $f$ is in the space $L^{p,q}_{w^*}(X,B)$, then we define
the norm $ \Vert f\Vert _{L^{p,q}(X,B^*)}$ to be the infimum of the
$L^{p,q}_{w^*}(X)$-norms of all nonnegative functions $g$ such that
$\Vert f\Vert _B \leq g$ almost everywhere on $X$. It can be proven that this
defines a Banach space structure (see Chapter $1$ in \cite{CeMe}). If
$f\dvtx X \rightarrow B$ is Bochner-measurable, and $B^*$ is separable, then
we say
that $f$ is in $L^{p,q}(X,B)$ if the \textit{measurable}
function
\[
x \mapsto\Vert f(x)\Vert _B
\]
is in $L^{p,q}(X)$. We will refer to $L^{p,q}(X,B)$ as the \textit
{Bochner--Lorentz space} with indices $p$ and $q$.

\subsection{Random semimetric spaces}

We will recall some basic notions from the ergodic theory of
subadditive cocycles. Classically, a \textit{subadditive cocycle} over a~measurable $\mathbb{Z}$-action $T$ on a
probability measure $(X,\mathfrak{F},\mu)$ is a measurable map $a \dvtx
\mathbb{Z}
\times X \rightarrow\mathbb{R}$ such that
\[
a(n+m,x) \leq a(n,x) + a(m,T_n x)\qquad \forall n,m \in\mathbb{Z}.
\]
A celebrated theorem of Kingman \cite{Ki} asserts that if $a(n,\cdot)$
is integrable with respect to $\mu$ for all $n$ in $\mathbb{Z}$, then there
exists a $T$-invariant real-valued measurable
function $A$ on $X$ such that
\[
\lim_{n \rightarrow+\infty} \frac{a(n,x) - n A(x)}{n} = 0
\]
almost everywhere on $(X,\mu)$. If the action $T$ is assumed to be
ergodic, then $A$ is necessarily constant. Furthermore, in this case,
\[
A = \inf_{n > 0} \frac{1}{n} \int_{X} a(n,x)\, d\mu(x).
\]
In this paper, we will be concerned with a generalization of this
theorem to measurable $\mathbb{Z}^d$-actions. We will need the
following definition.
\begin{defn}[(Random semimetric)]
Let $G$ be a locally compact and second countable group. Suppose that
$(X,\mathfrak{F})$ is a measurable space on which~$G$ acts measurably
and with an invariant probability
measure $\mu$. Let $(Y,\nu)$ be a $\sigma$-finite measure space, where
$\nu$ is a quasi-invariant measure under the action of $G$. A \textit
{random semimetric on $Y$, modeled
on the $G$-space $X$}, is a map $ \rho\dvtx X \times Y \times Y
\rightarrow
[0,\infty) $ such that the following conditions hold:
\begin{longlist}[(iii)]
\item[(i)] (symmetry) for all $x \in X$ and $y,y'$ in $Y$,
\[
\rho_x(y,y') = \rho_x(y',y) \quad \mbox{and}\quad  \rho_x(y,y) = 0;
\]
\item[(ii)] (triangle inequality) for all $x \in X$ and $y,y',y''$ in $Y$,
\[
\rho_x(y,y') \leq\rho_x(y,y'') + \rho_x(y'',y');
\]
\item[(iii)] (equivariance) for all $x \in X$ and $g \in G$ and $y,y'$
in $Y$,
\[
\rho_{gx}(y,y') = \rho_x(gy,gy').
\]
\end{longlist}
\end{defn}

\begin{rema*}
Let $(Z,d)$ be a metric space and suppose that $c \dvtx G \times X
\rightarrow
\operatorname{Isom}(Z,d)$ is a measurable map which satisfy the
equations
\[
c(gg',x) = c(g,x) c(g',gx)\qquad \forall g,g' \in G \mbox{ and } x \in X.
\]
It is easy to see that
\[
\rho_x(g,g') = d(c(g,x).z_0,c(g',x).z_0)
\]
defines a random semimetric on $G$, modeled on the $G$-space $X$, for
any choice of base point $z_0$ in $Z$. Indeed, by the cocycle property
of $c$, we have
\begin{eqnarray*}
\rho_{gx}(g',g'') &=& d(c(g',gx).z_0,c(g'',gx)) \\
&=& d(c(g,x)c(g',gx).z_0,c(g,x)c(g'',gx).z_0) \\
&=& d(c(gg',x).z_0,c(gg'',x)z_0) \\
&=& \rho_x(gg',gg'')
\end{eqnarray*}
for all $x \in X$ and $g,g',g''$ in $G$.
\end{rema*}

\subsection{Classical first passage percolation} \label{class}

First passage percolation was first defined by Hammersley and Welsh in
\cite{HaWe} and has served as one of the main inspirations for the
early developments of
subadditive ergodic theory. Let $(X,\mathfrak{F},\mu)$ be a probability
space on which the group $\mathbb{Z}^d$ acts ergodically and preserves the
measure $\mu$. We denote the
action by $T$. Let $f_1,\ldots,f_d$ be nonnegative measurable functions
on $X$ and define, for an edge $\bar{e} = (n,n+e_k)$ in the standard
Cayley graph of $\mathbb{Z}^d$, the weight
\[
t_x(\bar{e}) = f_k(T_{n}x), \qquad x \in X, k \in\{ 1,\ldots,d \},
\]
where $e_k$ denotes the $k$th standard basis vector in $\mathbb{Z}^d$. We
define the weight $t_x(\gamma)$ of a path $\gamma$ by summing the
individual weights on the edges of the path. For two points $m,n$ in
$\mathbb{Z}
^d$, we define
\[
\rho_x(m,n) = \inf\{ t_x(\gamma) \mid \mbox{$\gamma$ is a path from
$m$ to $n$} \}.
\]
It is clear from the construction that this defines a measurable map
from $X$ into the convex cone of semimetrics on $\mathbb{Z}^d$,
equipped with
the Borel structure coming from the topology of
pointwise convergence. Note that the relation $t_{T_k x}(\gamma) =
t_x(\gamma+ k)$ for $k \in\mathbb{Z}^d$ implies that
\[
\rho_{Tx}(m,n) = \rho_x(m+k,n+k)
\]
and thus $\rho$ is a random semimetric on $\mathbb{Z}^d$, modeled on
the $\mathbb{Z}
^d$-space $(X,\mu)$. By construction, the semimetric $\rho$ is inner.
The random semimetric space $(\mathbb{Z}^d,\rho)$
modeled on the $\mathbb{Z}^d$-space $X$ is known as the \textit
{classical first
passage percolation} model.

Note that in the case when $d=1$, we essentially recover the absolute
value of the \textit{Birkhoff sum}
\[
\sum_{k=0}^{n-1} f_1(T_k x)
\]
and thus the almost sure asymptotic behavior of the random semimetric
can be analyzed
using Birkhoff's ergodic theorem. When $d \geq2$, the situation is
more involved and new techniques are needed. The main part of this
paper is concerned with a generalization to general
semimetrics on $\mathbb{Z}^d$ of the following theorem of Boivin \cite{Bo90}.
\begin{thm}[(Boivin)] \label{Bo}
Suppose that $f_1,\ldots,f_d$ are in $L^{d,1}(X)$. There is then a
seminorm $L$ on $\mathbb{R}^d$ such that
\[
\lim_{n \rightarrow\infty} \frac{\rho_x(0,n)-L(n)}{|n|} = 0
\]
almost everywhere on $(X,\mu)$.
\end{thm}

\begin{rema*}
This theorem had previously been established for independent and
identically distributed edge-weights by Cox and Durrett \cite{CoDu}
($d=2$) and by Kesten~\cite{KeP} ($d \geq2$) under
weaker integrability conditions. However, it can be shown \cite{BoDe90}
that $L^{d,1}(X)$ is a sharp condition in the general ergodic case.
\end{rema*}

The definition of classical first passage percolation described above
extends naturally to a more general situation. Let $G$ be a finitely
generated group and suppose that $S$ is
a finite subset of $G$ such that $S$ and $S^{-1}$ are disjoint and $S
\cup S^{-1}$ generates~$G$ as a group. Suppose that $ \{ f_s \}_{s \in
S}$ is a set of nonnegative measurable
functions on a probability measure space $(X,\mu)$ with a
measure-preserving \textit{right} action by $G$. For every $g$ in $G$ and
edge $(g,gs)$ in the Cayley graph of
$(G,S\cup S^{-1})$, we define the random weight $t_x(g,gs) = f_s(xg)$.
In analogy with the scheme above, we define the distance $\rho$ between
two points $g$ and $g'$ in $G$
to be the infimum of the weights over all paths between $g$ and $g'$.
By construction, $\rho$ is a semimetric and
\[
\rho_x(hg,hg') = \rho_{xh}(g,g')
\]
for all $g,g',h$ in $G$ and $x$ in $X$. It is not clear that Boivin's
\textit{proof} of Theorem \ref{Bo} immediately extends to the case when
$G = \mathbb{Z}^d$ and $S$ is \textit{not} the
standard generating set. Note, however, that Theorem \ref{main2} covers
this case.

\subsection{Cohomology of Borel groupoids} \label{coho}

In this subsection, we will define various important types of cocycles.
A more conceptual explanation can be given in the language of
groupoids; however, we will refrain from making
very general statements and will restrict our attention to the first
order cohomology of a groupoid.

\begin{defn}[(Borel cocycle)]
Let $(Z,d)$ be a metric space and $G$ a topological group. Suppose that
$X$ is a $G$-space. A map $c\dvtx G \times X \rightarrow\operatorname
{Isom}(Z,d)$ such that
\[
(g,x) \mapsto c(g,x).z
\]
is measurable for all $z$ in $Z$, with respect to the Borel $\sigma
$-algebra on $Z$, and
\[
c(gg',x) = c(g,x) c(g',gx) \qquad \forall g,g' \in G, x \in X,
\]
is called a \textit{Borel cocycle} over the $G$-space $X$.
\end{defn}

\begin{defn}[(Gelfand cocycle)]
Let $B$ be a Banach space and $G$ a locally compact and second
countable group. Let $(X,\mu)$ be a probability measure space, where
$\mu$ is a $G$-invariant measure
and $(Y,\nu)$ a $\sigma$-finite measure space, where $\nu$ is a
quasi-invariant measure under the action of $G$. Let $c \dvtx G \times X
\rightarrow\operatorname{Isom}(B^{*})$ be a Borel cocycle.
A map $s \dvtx X \times Y \times Y \rightarrow B^*$ is called a \textit{Gelfand
$B$-cocycle with respect to the Borel cocycle $c$} if the following
conditions hold:
\begin{longlist}[(iii)]
\item[(i)] (additivity) for all $x \in X$,
\[
s_x(y,y'') + s_x(y'',y') = s_x(y,y')\qquad \forall y,y',y'' \in Y;
\]
\item[(ii)] (equivariance) for all $x \in X$,
\[
c(g,x).s_{gx}(y,y') = s_x(gy,gy')\qquad \forall y,y' \in Y
\mbox{ and } g \in G;
\]
\item[(iii)] (measurability) the maps
\[
x \mapsto s_x(y,y') \quad \mbox{and}\quad  x \mapsto\Vert s_x(y,y')\Vert _{B^*}
\]
are weak*-measurable for every $y,y' \in Y$.
\end{longlist}
\end{defn}

We say that $s$ is in the Lorentz space $L^{p,q}_{w*}(X,B)$ if the map
$x \mapsto\Vert s_x(y,y')\Vert _{p,q}$ is in $L^{p,q}(X)$ for all $y,y'$ in $Y$.

\begin{rema*}
If the cocycle is trivial, that is, if there is an isometric
representation $\pi$ of $G$ on $B^{*}$ such that $c(g,x) = \pi(g) $ for
all $x \in X$ and $g \in G$,
we will refer to $s$ as a \textit{Gelfand $B$-cocycle} with respect to
the representation $\pi$.
\end{rema*}

We also define two related types of cocycles, where stronger versions
of measurability are assumed.
\begin{defn}[(Pettis cocycle)]
A map $s \dvtx X \times\mathbb{Z}^d \times\mathbb{Z}^d \rightarrow B$
is called a~\textit{Pettis
$B$-cocycle with respect to the Borel cocycle $c$} if it is a Gelfand
cocycle with respect to the cocycle $c$ and the maps
\[
x \mapsto s_x(y,y')
\]
are weakly measurable for all $y,y'$ in $Y$.
\end{defn}
\begin{rema*}
Note that, in the definition of a Pettis cocycle, we do not insist that
$s$ takes values in the \textit{dual} of a Banach space $B$. Thus, the
formulation of the definition is slightly misleading, but we hope that
this will not cause any confusion for the reader.
\end{rema*}

In Section \ref{bochner}, we will need the following cocycles, which
are measurable in a strong sense.
\begin{defn}[Bochner cocycle]
A map $s \dvtx X \times\mathbb{Z}^d \times\mathbb{Z}^d \rightarrow B$
is called a~\textit{Bochner
$B$-cocycle with respect to the Borel cocycle $c$} if it is a Gelfand
cocycle with respect to $c$ and the maps
\[
x \mapsto s_x(y,y')
\]
are Bochner measurable for all $y,y'$ in $Y$.
\end{defn}

\begin{rema*}
If $B$ is a separable Banach space, it follows from Pettis's
measurability theorem (see, e.g., Chapter $1$ of \cite{Di}) that every
Pettis cocycle is
a Bochner cocycle. The converse is obvious.
\end{rema*}

One connection between Gelfand $B$-cocycles and random semimetrics on
$Y$ is suggested by the following proposition.
\begin{pro}
Let $G$ be a locally compact group. Suppose that $s \dvtx X \times Y \times
Y \rightarrow B^{*}$ is a Gelfand $B$-cocycle with respect to a Borel cocycle
$c$. Then
\[
\rho_x(y,y') = \Vert s_x(y,y')\Vert _{B^*},\qquad  y,y' \in Y,
\]
is a random semimetric on $Y$, modeled on the $G$-space $X$.
\end{pro}

\begin{pf}
The measurability is clear from the definition of $s$. From the
additivity property of $s$, it follows that
\[
s_x(y,y') = -s_x(y',y) \qquad \mbox{for all $x,y,y'$.}
\]
Thus, $\rho_x$ is symmetric and $\rho_x(y,y) = 0$.
For the triangle inequality, we observe that
\begin{eqnarray*}
\Vert s_x(y,y')\Vert _{B^{*}} &=& \Vert s_x(y,y'') + s_x(y'',y')\Vert _{B^{*}} \\
&\leq& \Vert s_x(y,y'')\Vert _{B^{*}} + \Vert s_x(y'',y')\Vert _{B^{*}} \\
&=& \rho_x(y,y'') + \rho_x(y'',y')
\end{eqnarray*}
for all $y,y',y''$ in $Y$. Finally, to prove equivariance, we note that
since $c$ takes values in the isometry group of $B^{*}$,
\begin{eqnarray*}
\rho_{gx}(y,y') &=& \Vert s_{gx}(y,y')\Vert _{B^{*}} =
\Vert c(g,x).s_{gx}(y,y')\Vert _{B^{*}} \\
&=& \Vert s_{x}(gy,gy')\Vert _{B^*} = \rho_x(gy,gy')
\end{eqnarray*}
for all $g \in G$ and $y,y'$ in $Y$. In the second-to-last equality,
the equivariance property of $s$ was used.
\end{pf}

\begin{rema*}
We will see in Theorem \ref{thm:structure} that these examples of
random semimetrics are the only such examples.
This observation will be one of the crucial steps in the proof of
Theorem \ref{thm:main}.
\end{rema*}

\subsection{Representation of subadditive cocycles}

In this subsection, we will prove the following, important, structure theorem.

\begin{thm} \label{thm:structure}
Let $G$ be a locally compact, second countable group. Suppose that
$(X,\mu)$ is a probability measure space
with a $G$-invariant ergodic measure $\mu$. Suppose that $(Y,\nu)$ is a
$G$-space with a quasi-invariant $\sigma$-finite
measure. If $\rho$ is a random $G$-equivariant semimetric on $Y$,
modeled on the $G$-space $(X,\mu)$, then there exists a
Gelfand $L^{1}(Y,\nu)$-cocycle, with respect to the left-regular
representation $\lambda$ of $G$ on $L^{\infty}(Y,\nu)$, on the $G$-space
$(X,\mu)$ such that
\[
\rho_x(y,y') = \Vert s_x(y,y')\Vert _{L^{\infty}(Y,\nu)}.
\]
\end{thm}

\begin{pf}
The proof is based on the following trivial observation:
\[
\rho_x(y,y') = \sup_{y'' \in Y} | \rho_x(y,y'') - \rho_x(y'',y') |,
\]
which is a direct consequence of the triangle inequality. We define
\[
s_x(y,y') = \rho_x(y,\cdot) - \rho_x(\cdot,y') \in L^{\infty
}(Y,\nu).
\]
Note that
\[
s_x(y,y'') + s_x(y'',y') = \rho_x(y,\cdot) - \rho_x(\cdot,y'') +
\rho
_x(\cdot,y'') - \rho_x(\cdot,y') = s_x(y,y')
\]
and that
\[
\lambda(g).s_{gx}(y,y') = \rho_{gx}(y,g^{-1}\cdot) - \rho
_{gx}(g^{-1}\cdot,y') = s_x(gy,gy').
\]
To prove measurability, we first note that the map $x \mapsto
\Vert s_x(y,y')\Vert _{L^{\infty}(Y,\nu)}$ is measurable,
by definition. Thus, we only need to prove that the map $s_x(y,y')$ is
weak*-measurable. If we choose
$\eta\in L^{1}(Y,\nu)$, then
\[
\langle\eta, s_x(y,y') \rangle= \int_{Y} \bigl(\rho_x(y,y'')-\rho
_x(y'',y')\bigr) \eta
(y'')\, d\nu(y'')
\]
is measurable by Fubini's theorem, since, by definition, the map
\[
(x,y'') \mapsto\bigl(\rho_x(y,y'')-\rho_x(y'',y')\bigr) \eta(y'')
\]
is clearly measurable on the probability measure space $ ( X \times Y ,
\mu\times\nu) $ for almost every choice of $y,y' \in Y$
with respect to $\nu\times\nu$.
\end{pf}

\begin{rema*}
In the paper \cite{Ki}, Kingman asked the naturally arising question as
to whether every subadditive cocycle $a$ on $\mathbb{Z}$-space $X$ has
a representation
of the form
\[
a(n,x) = \sup_{i \in I} \sum_{k=0}^{n-1} f_i(T^k x),\qquad n \in\mathbb{N},
\]
where $ \{ f_i \}_{i \in I}$ is a set of real-valued measurable
functions on $X$ and $I$ is some countable index set. This question was
later answered in the negative by
Hammersley in \cite{Ha}. Theorem \ref{thm:structure} gives a positive
answer to an extended version of Kingman's question, where the
functions $f_i$ are
allowed to be Banach-space-valued and the supremum is replaced by the
corresponding Banach norm. However, it is certainly an inconvenience
that the proof
requires the Banach space $B^*$ to be nonseparable. Therefore, it seems
appropriate to ask the following question.
\begin{qu*}
Can every random $G$-equivariant semimetric $\rho$ on a $G$-space $Y$,
quasi-invariant under the action of $G$ and modeled on a $G$-space
$(X,\mu)$, be represented as the norm
of a Gelfand $B$-cocycle, where
$B^*$ is a \textit{separable} Banach space?
\end{qu*}
\end{rema*}

\subsection{Asymptotic shape theorems}

We will now outline the main steps in the proof of Theorem \ref{main2}.
We first make some preliminary observations and remarks.
\begin{pro}
Suppose that $\rho$ is a random $G$-equivariant semimetric on~$Y$,
modeled on the $G$-space $(X,\mu)$. The function
\[
r(y,y') = \int_{X} \rho_x(y,y')\, d\mu(x)
\]
is then a $G$-invariant semimetric on $Y$.
\end{pro}

\begin{pf}
The axioms for a semimetric are easily verified. The rest of the proof
consists of the following simple
calculation:
\[
r(gy,gy') = \int_{X} \rho_x(gy,gy')\, d\mu(x) = \int_{X} \rho
_{gx}(y,y')\, d\mu(x) = r(y,y').
\]
\upqed\end{pf}

The study of the almost sure asymptotic geometry of random semimetric
spaces will be referred to as \textit{generalized first passage
percolation}. We begin by
describing some general features of this theory. Suppose that $Y$ is a
locally compact space and that $r$ is dominated by a $G$-invariant
metric $\eta$ such that
\[
\liminf_{y \rightarrow\infty} \eta(y,o) = + \infty
\]
for every choice of base point $o \in Y$. We say that the random
$G$-equivariant semimetric on $Y$ satisfies an \textit{asymptotic
shape theorem} (with respect to the $G$-invariant metric $\eta$) if
there exists a measurable function $L \dvtx Y \rightarrow[0,\infty) $
such that
\[
\limsup_{y \rightarrow\infty} \biggl| \frac{\rho_x(o,y)-L(y)}{\eta
(o,y)} \biggr|
= 0.
\]
This paper is concerned with a general asymptotic shape theorem for
actions of the group $\mathbb{Z}^d$ on probability spaces. We will
specialize the
above situation to the case where $G = \mathbb{Z}^d$ and $Y = \mathbb
{Z}^d$. In this
case, $L$ can be realized as a norm on $\mathbb{R}^d$ and $\eta$ will
be taken to be the standard word-metric on $\mathbb{Z}^d$.

The following important lemma is proved in \cite{Bo90}.
\begin{lem}[(Boivin's lemma)]
Suppose that $\rho$ is a random $\mathbb{Z}^d$-equivariant semimetric
on $\mathbb{Z}
^d$, modeled on an ergodic $\mathbb{Z}^d$-space $(X,\mu)$, where $\mu
$ is a
probability measure. If there is a positive constant $C$ such that
\[
\mu\biggl( \biggl\{ x \in X \mbox{ }\Big|\mbox{ } \sup_{n \neq0} \frac{\rho_x(0,n)}{|n|} \geq
\lambda \biggr\} \biggr) \leq\frac{C}{\lambda^d}\qquad  \forall\lambda> 1,
\]
then there exists a seminorm $L$ on $\mathbb{R}^d$ such that
\[
\lim_{|n| \rightarrow\infty} \frac{\rho_x(0,n) - L(n)}{|n|} = 0
\]
almost everywhere on $(X,\mu)$.
\end{lem}

Thus, in order to prove an asymptotic shape theorem, we will need a
maximal inequality. Let $s$ be a Gelfand $B$-cocycle and
define the function
\[
Ms(x) = \sup_{n \neq0} \frac{\Vert s_x(0,n)\Vert _{B^*}}{|n|},\qquad x \in X.
\]

 We prove the following maximal inequality.
\begin{thm}[(Maximal inequality)] \label{maxineqthm}
Let $B$ be a Banach space. Suppose that~$s$ is a Gelfand $B$-cocycle on
an ergodic $\mathbb{Z}^d$-space $(X,\mu)$, where $\mu$ is a probability
measure. Suppose that $s(0,n)$ is in $ L^{d,1}_{w^{*}}(X,\mu,B^*)$ for
every $n \in\mathbb{Z}^d$. There exists a positive constant $C$ such that
\[
\mu\bigl( \{ x \in X \mid Ms(x) \geq\lambda \} \bigr) \leq
\frac{C}{\lambda^{d}} \Vert s\Vert _{L^{d,1}_{w^*}(X,B)}
\]
for all $ \lambda\geq1 $.
\end{thm}

 The proof of this theorem will be presented in Section \ref{maxineq}.
An immediate corollary of this result is
the following theorem.
\begin{thm}[(Asymptotic shape theorem)] \label{thm:main}
Suppose that $\rho$ is a random $\mathbb{Z}^d$-semimetric modeled on
an ergodic
$\mathbb{Z}^d$-space $(X,\mu)$, where $\mu$ is a probability measure.
Suppose that $\rho(0,n)$ is in $L^{d,1}(X,\mu)$ for every $n \in
\mathbb{Z}^d$.
There then exists a~seminorm $L$ on $\mathbb{R}^d$ such that
\[
\lim_{|n| \rightarrow\infty} \frac{\rho_x(0,n) - L(n)}{|n|} = 0
\]
almost everywhere on $(X,\mu)$.
\end{thm}

\begin{rema*}
This theorem was proven by Boivin in \cite{Bo90} in the case of certain
\textit{inner} random semimetrics on $\mathbb{Z}^d$. The proof is slightly
different and does not
seem to extend to the general situation. Note that when $d=1$, Boivin's
theorem is essentially equivalent to Birkhoff's ergodic theorem, while
our theorem is
strictly stronger.
\end{rema*}

\subsection{Maximal inequalities} \label{maxineq}

The goal of this section is to establish Theorem~\ref{maxineqthm}. The
proof closely follows the arguments outlined by Derriennic
and Boivin in \cite{BoDe90}. We begin to recall the basic combinatorial
lemma used by Boivin and Derriennic in their proof. A detailed
proof can be found in \cite{BoDe90}.
\begin{lem} \label{comb}
For every $n \in\mathbb{Z}^d$, $n \neq0$, let
\[
P_n = \{ m \in\mathbb{Z}^d \mid |n-m| \leq|n|/2 \}.
\]
Let $H$\vspace*{-2pt} be a coordinate hyperplane of $\mathbb{Z}^d$ such that $H \cap
P_x =
\varnothing$. Let
$(H_j)_{j=1}^{d-1}$ be an increasing sequence of coordinate subspaces
of $\mathbb{Z}^d$ such
that $\dim{H_j} = j$ and $H_{d-1} = H$. There exists a set $\mathcaligr
{E}_n$ of elementary
paths $\gamma$ in $\mathbb{Z}^d$, joining $0$ to $n$, such that:
\begin{longlist}[(iii)]
\item[(i)] the cardinality of $ \mathcaligr{E}_n $ is $|n|^{d-1}$;
\item[(ii)] each $\gamma$ is entirely included in the set
\[
\{ m \in\mathbb{Z}^d \mid |m| \leq2|n| \};
\]
\item[(iii)] for every $m \in P_n$ and $m \neq n$,
\[
\bigl| \{ \gamma\in\mathcaligr{E}_n | m \in\gamma \} \bigr| \leq C\biggl(\frac
{|n|}{|n-m|}\biggr)^{d-1};
\]
\item[(iv)] for every $m \notin P_n $ with $|m| \leq2|n|$,
\[
\bigl| \{ \gamma\in\mathcaligr{E}_n | m \in\gamma \} \bigr| \leq|n|^{d-j(m)}
\]
\item[]for $j(m) = \sup\{ j = 1 , \ldots, d-1 \mid m \in H_j \} $; \\
\item[(v)] for every $m \notin H \cup P_n $ with $|m| \leq2|n|, $
\[
\bigl|\{ \gamma\in\mathcaligr{E}_x | m \in\gamma \}\bigr| \leq1.
\]
\end{longlist}
\end{lem}

For an elementary path $\gamma= \{ n_1,\ldots,n_r \}$ between $0$ and
$n$ in $\mathbb{Z}^d$, we define
\[
A_\gamma f(x) = \sum_{k=0}^{r-1} f(T_{n_k}x),
\]
where $f \dvtx X \rightarrow\mathbb{R}$ is a measurable function. The
following lemma was
proven in~\cite{BoDe90}.
\begin{lem} \label{help}
Suppose that $f$ is a nonnegative and measurable function on $X$. Then,
for all $n \neq0$,
\begin{eqnarray*}
\frac{1}{|\mathcaligr{E}_n|} \biggl| \sum_{\gamma\in\mathcaligr{E}_n}
A_\gamma f(x)\biggr | &\leq& C
\biggl[ \sum_{H \in\mathcaligr{H}} \frac{1}{|n|^{\dim{H}}} \mathop{\mathop{\sum}_{
{m \in H}}}_{ |m| \leq2|n|} f(T_m x)\\
&&\quad{} + \frac{1}{|n|} \biggl(f(T_nx) + \sum_{0 < |n-m| \leq|n|/2} \frac{
f(T_m x)}{|m-n|^{d-1}}\biggr) \biggr],
\end{eqnarray*}
where $\mathcaligr{H}$ denotes the collection of all coordinate subspaces
of $\mathbb{Z}^d$.
\end{lem}

This lemma readily implies the following estimate.
\begin{pro} \label{mainest}
For every nonzero $n \in\mathbb{Z}^{d}$, we have
\begin{eqnarray*}
\frac{\Vert s_x(0,n)\Vert _{B^*}}{|n|} & \leq& C \biggl[ \sum_{H \in\mathcaligr
{H}} \frac{1}{|n|^{\dim{H}}} \mathop{\mathop{\sum}_{m \in H}}_{|m| \leq2|n|}
f(T_m x) \\
&&\quad{}+ \frac{1}{|n|} \biggl(f(T_nx) + \sum_{0 < |n-m| \leq|n|/2} \frac{
f(T_m x)}{|m-n|^{d-1}}\biggr) \biggr],
\end{eqnarray*}
where $\mathcaligr{H}$ denotes the collection of all coordinate subspaces
of $\mathbb{Z}^d$ and
\[
f(x) = \sup_{k=1,\ldots,d} \max{ (\Vert s_x(0,e_{k})\Vert _{B^{*}},
\Vert s_x(0,-e_{k})\Vert _{B^{*}})}.
\]
\end{pro}

\begin{pf}
For every elementary path $\gamma_n = \{ n_1,\ldots,n_{r} \}$ from $0$
to $n$, we write
\[
s_x(0,n) = \sum_{k=0}^{r-1} s_x(n_k,n_{k+1}) = \sum_{k=0}^{r-1}
\lambda
(n_k).s_{T_{n_k}x}(0,n_{k+1}-n_k),
\]
where $\lambda$ is the left-regular representation of $L^{\infty
}(\mathbb{Z}
^d)$. Thus, since $|n_{k+1}-n_k| = 1$ for all $k$,
we have
\[
\Vert s_x(0,n)\Vert _{B^*} \leq\sum_{k=0}^{r-1} f(T_{n_k}x) = A_\gamma f(x).
\]
We now take the average over the set $\mathcaligr{E}_n$. By Lemma \ref
{comb} and Proposition \ref{help}, we have
\begin{eqnarray*}
\frac{\Vert s_x(0,n)\Vert _{B^*}}{|n|} &\leq& C \biggl[ \sum_{H \in\mathcaligr{H}}
\frac{1}{|n|^{\dim{H}}} \mathop{\mathop{\sum}_{m \in H}}_{|m| \leq2|n|}
f(T_m x) \\
&&\quad{}+ \frac{1}{|n|} \biggl(f(T_nx) + \sum_{0 < |n-m| \leq|n|/2} \frac{
f(T_m x)}{|m-n|^{d-1}}\biggr) \biggr].
\end{eqnarray*}
%
\upqed\end{pf}

The following lemma was proven in \cite{BoDe90} for general actions of
$\mathbb{Z}^d$ on probability spaces.
\begin{lem} \label{mainest2}
Suppose that $f$ is a nonnegative and measurable function on $X$. There
is then a constant $C > 0$ such that
\[
\mu\biggl( x \in X \Big| \sup_{n \in\mathbb{Z}^d} \frac{1}{|n|}\biggl( f(T_n x)
+ \sum_{0 < |n-m| \leq|n|/2} \frac{f(T_m x )}{|n-m|^{d-1}} > \lambda
\biggr) \biggr)
\leq C \biggl(\frac{1}{\lambda} \Vert f\Vert _{d,1}\biggr)^{d}
\]
for all $\lambda\geq1$.
\end{lem}

 Define the maximal function
\[
Ms(x) = \sup_{n \neq0} \frac{\Vert s_x(0,n)\Vert _{B^*}}{|n|}
\]
for a Gelfand $B$-cocycle $s$. The maximal inequality for the first
terms in the estimate in Lemma \ref{help} are taken care of by Wiener's
maximal inequality \cite{Wi}.
Proposition \ref{mainest} and Lemma \ref{mainest2} now imply the
following theorem.
\begin{thm}[(Maximal inequality)]
Let $B$ be a Banach space. Suppose that~$s$ is a Gelfand $B$-cocycle on
an ergodic $\mathbb{Z}^d$-space $(X,\mu)$, where $\mu$ is a probability
measure. Suppose that $s(0,n)$ is in $ L^{d,1}_{w^{*}}(X,\mu,B^*)$ for
every $n \in\mathbb{Z}^d$. There then exists a positive constant $C$
such that
\[
\mu\bigl( \{ x \in X | Ms(x) \geq\lambda \} \bigr) \leq
\biggl( \frac{C}{\lambda} \Vert s\Vert _{L^{d,1}_{w^*}(X,B)} \biggr)^d
\]
for all $ \lambda\geq1 $.
\end{thm}

\subsection{Ergodic theorems for Bochner cocycles} \label{bochner}

In this subsection, we will be concerned with a slight generalization
of the ergodic theorem of Boivin and Derriennic in \cite{BoDe90} to
vector-valued
cocycles. We will see an application of this theorem to horofunctions
in random media in Section \ref{horo}. The main ingredient of the proof
is a
result of Phillips, ensuring that the Bochner--Lorentz space
$L^{d,q}(X,B)$ is a reflexive Banach space if $ 1 < q < \infty$. This
will allow for standard
splitting theorems to be used (see Chapter $2$ of Krengel's book \cite
{Kr} for more references).

\begin{thm} \label{thm:linear}
Let $B$ be a reflexive Banach space. Suppose that $s$ is Bochner
$B$-cocycle on an ergodic $\mathbb{Z}^d$-space $(X,\mu)$, where $\mu
$ is a
probability measure. Let $q > 1$ and suppose that the cocycle $s(0,n)$
is in $L^{d,q}_s(X,B^*)$ for every $ n \in\mathbb{Z}^d$. There is
then a
linear and continuous
map $L\dvtx \mathbb{R}^d \rightarrow B^* $ such that
\[
\lim_{|n| \rightarrow\infty} \frac{s_x(0,n)-L(n)}{|n|} = 0
\]
almost everywhere on $(X,\mu)$.
\end{thm}

Before we turn to the proof of Theorem \ref{thm:linear}, we recall some
basic splitting theorems for $\mathbb{Z}^d$-actions on
Bochner--Lorentz spaces. The following theorem is due to Phillips \cite{Phi43}.
\begin{thm}
Let $(X,\mathcaligr{F},\mu)$ be a measure space. Suppose that $B$ is
a reflexive
Banach space, $ 1 \leq p < \infty$ and $ 1 < q < \infty$.
Then $L^{p,q}(X,B)$ is reflexive.
\end{thm}

 By the well-developed splitting theory of semigroups of isometries on
reflexive Banach spaces (see, e.g., Chapter $2$ of \cite{Kr}), this
implies that every Bochner cocycle
$s$ in $L^{d,1}(X,B^*)$ can be written as a limit in $L^{d,1}(X,B)$ of
\[
s = \lim_{j \rightarrow\infty} r + c^{j},
\]
where $r$ is an invariant cocycle and $c^{j}$ is a sequence of
coboundaries, that is, cocycles of the form
\[
c^{j}_x(0,n) = g_j(x) - g_j(T_n x),\qquad g_j \in L^{d,q}(X,B^*), q
> 1,
\]
and extended by equivariance.

\begin{pf*}{Proof of Theorem \protect\ref{thm:linear}}
Note that the theorem is trivial for invariant cocycles and
coboundaries. By Banach's principle and Theorem \ref{maxineqthm},
the set of all cocycles for which the theorem holds is closed in
$L^{d,q}(X,B)$. Since the span of invariant cocycles and
coboundaries is dense in $L^{d,q}(X,B^*)$, we are done.
\end{pf*}

\section{Applications}

\subsection{Random Schr\"odinger operators}

Let $(X,\mathcaligr{F},\mu,T)$ be an ergodic $\mathbb{Z}$-space and
suppose that $S$ is a
measurable map from $\mathbb{Z}\times X \rightarrow GL_d(\mathbb{R})$
which satisfies the equations
$S(0,\cdot) = I$ and
\[
S(n+m,x) = S(n,T^m x) S(m,x) \qquad\forall m,n \in\mathbb{Z}
\]
almost everywhere on $X$. The asymptotic behavior of the random semimetric
\begin{eqnarray}
\rho_x(m,n) = \max{\bigl(\log^{+}(\Vert S_n(x) S_m(x)^{-1}\Vert) ,\log
^{+}(\Vert S_m(x)S_n(x)^{-1}\Vert )\bigr)}\nonumber\\
\eqntext{\forall n,m \in\mathbb{Z},}
\end{eqnarray}
has been the subject of a detailed study in the theory of random Schr\"
odinger operators over the years. The first convergence result, prior
to Kingman's paper, was due to Furstenberg and Kesten \cite{KF}, where
the almost sure limit
\[
A = \lim_{n \rightarrow\infty} \frac{\rho_x(0,n)}{|n|}
\]
was established. Before we begin our discussion of multiparameter
analogs, we first describe the connection to random Schr\"odinger operators.
Let\break  $(X,\mathcaligr{F},\mu,T)$ be an ergodic
probability-measure-preserving system
and suppose that $V$ is a real-valued measurable function on $X$.
We consider, for a fixed $x$ in $X$, the following discrete analog of
the Schr\"odinger equation:
\[
v_{n+1} + v_{n-1} + V(T^n x) v_n = \lambda v_n\qquad \forall n \in
\mathbb{Z},
\]
with $v_o=a$ and $v_1=b$, where $\lambda$ is assumed to be real. If we
introduce the vectors $u_n = (v_n , v_{n+1})^t$, the equation can be
written in the equivalent form
\[
u_{n+1} =
\pmatrix{
0 & 1 \cr
-1 & \lambda- V(T^n x) }
u_n,
\qquad
u_0 =
\pmatrix{
a \cr
b }
\]
and thus
\[
u_n = S(n,x) u_o\qquad \forall n \in\mathbb{Z},
\]
with $S$ is generated by
\[
S(1,x) =
\pmatrix{
0 & 1 \cr
-1 & \lambda- V(x)}.
\]
Hence, the generic [in terms of the measure space $(X,\mathcaligr
{F},\mu)$]
asymptotic behavior of the solutions of random Schr\"odinger operators
on $\mathbb{Z}$ is governed by
the random semimetric $\rho_x$ defined above. By a remarkable tour de
force, Furstenberg and Kesten established the almost sure limit
\[
A = \lim_{n \rightarrow\infty} \frac{\rho_x(0,n)}{|n|}.
\]
This result predates Kingman's subadditive ergodic theorem and the
methods of Furstenberg and Kesten were indeed quite different from the
ones Kingman later used.

This example leads to a natural generalization for $\mathbb
{Z}^d$-actions. Let
$(X,\mathcaligr{F},\mu,T)$ be an ergodic $\mathbb{Z}^d$-space and
let $S \dvtx \mathbb{Z}^d \times X
\rightarrow GL_k(\mathbb{R})$ be
a measurable map which satisfies $S(0,\cdot) = I$ and
\[
S(n+m,x) = S(n,T^m x) S(m,x)\qquad  \forall m,n \in\mathbb{Z}^d,
\]
almost everywhere on $X$. Define the random semimetric
\begin{eqnarray}
\rho_x(m,n) = \max{\bigl(\log^{+}(\Vert S_n(x) S_m(x)^{-1} \Vert) ,\log
^{+}(\Vert S_m(x)S_n(x)^{-1}\Vert )\bigr)}\nonumber\\
\eqntext{\forall n,m \in\mathbb{Z}^d.}
\end{eqnarray}
Note that the sequence $u_n = S(T^n x) u$ is the solution of the random
difference equation
\[
\sum_{|e|=1} u_{n+e} = \biggl( \sum_{|e|=1} S(e,T^n x) \biggr) u_n
\qquad\forall n \in\mathbb{Z}^d,
\]
with $u_o = u$, where $|\cdot|$ denotes the $\ell^\infty$-metric on
$\mathbb{Z}
^d$. Note that the existence of a map $S$ with the above properties is
not obvious and, indeed, we do not expect any new examples for $d \leq
2$. However, any embedding of $\mathbb{Z}^d$ into $GL_k(\mathbb{R})$
for sufficiently
large $k$ will give rise to maps $S$ with the above properties and
hence the class of new difference equations which can be solved by this
method is nontrivial. In this class, the following theorem can be
deduced from Theorem \ref{main2}.
\begin{thm}
Suppose that for some $\varepsilon> 0$,
\[
\int_{X} ( \log^+\Vert S_n(x)\Vert  )^{d+\varepsilon}\, d\mu(x) < + \infty
\]
for all $|n| = 1$. There is then a seminorm $L$ on $\mathbb{R}^d$ with the
property that
\[
\lim_{n \rightarrow\infty} \frac{\rho_x(0,n)-L(n)}{|n|} = 0
\]
almost everywhere on $X$.
\end{thm}

The class of difference equations which can be solved by means of the
above scheme can probably be considerably enlarged if Theorem \ref
{main2} is extended to more
general linear groups, which motivates the study of extensions of
Boivin and Derriennic's result to more general groups.

\subsection{A multiplicative ergodic theorem}

In this subsection, we will establish a~multiplicative ergodic theorem
for general Pettis $\mathbb{Z}$-cocycles on ergodic probability
measure spaces
$(X,\mu)$ with values in separable Banach spaces
$B$ with the Radon--Nikodym property. The formulation is close to the
celebrated Karlsson--Ledrappier ergodic theorem \cite{KaLe}. However,
their paper is concerned with a special kind of
a Pettis cocycle with values in the Banach space of continuous
functions on an infinite compact metrizable space, which,
unfortunately, does not possess the Radon--Nikodym
property \cite{Di}. It is not unlikely that our ergodic theorem holds
in greater generality (e.g., nonseparable or weakly compact generated
Banach spaces). For the present proof and methods,
the Radon--Nikodym assumption seems to be sharp.

We begin by recalling the definition and some basic facts about Banach
spaces with the Radon--Nikodym property. Recall that a vector measure
$\nu$ is $\mu$\textit{-continuous} if
\[
\lim_{\mu(E) \rightarrow0} \nu(E) = 0.
\]

\begin{defn}[(Radon--Nikodym property)]
A Banach space $B$ has the \textit{Radon--Nikodym property with respect
to the measure space $(X,\mathcaligr{F},\mu)$} if, for each $\mu
$-continuous vector
measure $\nu\dvtx \mathcaligr{F}\rightarrow B$ of bounded variation, there
exists $g \in L^1(X,B)$ such that
\[
\nu(E) = \int_{E} g\, d\mu\qquad \forall E \in\mathcaligr{F},
\]
in the sense of Bochner integrals. A Banach space $B$ has the \textit
{Radon--Nikodym property} if $B$ has the Radon--Nikodym property with
respect to any finite measure space.
\end{defn}

 Classical examples of Banach spaces with the Radon--Nikodym property
include reflexive Banach space and Banach spaces with separable dual
spaces. Examples of Banach space without the Radon--Nikodym
property are $L^1([0,1])$ and $C(H)$, where $H$ is an infinite compact
Hausdorff space. The notion of a Radon--Nikodym space is now fairly
well understood and a very readable account of results and
techniques can be found in \cite{Di}.

We will need the following theorem by Bochner and Taylor \cite{BoPh}.
\begin{thm}
Let $(X,\mathcaligr{F},\mu)$ be a finite measure space, $B$ be a
Banach space and
$1\leq p < \infty$. Then $L^{p}(X,B)^* = L^{q}(X,B^*)$, where $\frac
{1}{p} + \frac{1}{q} = 1$, if and only if~$B^*$ has the
Radon--Nikodym property with respect to $\mu$.
\end{thm}

 In particular, this implies that if $B$ has the Radon--Nikodym
property, then
\[
\Vert f\Vert _{L^1(X,B)} = \sup_{\Vert \lambda\Vert _{\infty,B^*} \leq1} \int_{X}
\langle
\lambda(x) , f(x) \rangle\, d\mu(x)
\]
for every Bochner measurable function $f\dvtx X \rightarrow B $, where $\Vert
\cdot
\Vert _{\infty,B^*}$ denotes the $L^{\infty}(X,B^*)$-norm. More generally,
we will write $\Vert \cdot\Vert _{q,\mathcaligr{C}}$ if we restrict the
elements in $L^{q}(X,B^*)$ to take values in $\mathcaligr{C}$, for $q >
1$.

Suppose that $s$ is a Pettis $B$-cocycle on a probability measure space
$(X,\mu)$ with respect to a $\mathbb{Z}$-action $T$ and Borel cocycle $c$,
where the Banach space $B$ is supposed to be
separable and to have the Radon--Nikodym property. Note that this
implies that $s$ is also Bochner integrable. We also assume that the
function $x \mapsto\Vert s_x(0,n)\Vert _{B}$ is integrable for all $n \in
\mathbb{Z}$.
Suppose that there exists a weak*-compact subset $\mathcaligr{C}$ of
$B^{*}_1$ which is invariant under the dual action of the cocycle $c$,
such that
\[
\Vert s_x(m,n)\Vert _{1,B} = \sup_{\Vert \lambda\Vert _{\infty,\mathcaligr{C}} \leq1}
\langle\lambda, s_x(m,n) \rangle\qquad \forall m,n \in\mathbb{Z}.
\]

It is a well-known fact (see, e.g., Chapter V.5.1 of \cite{DS}) that
$\mathcaligr{C}$ is metrizable and thus separable. By subadditivity, the
following nonnegative limit exists:
\[
A := \lim_{n \rightarrow\infty} \frac{1}{n} \Vert s(0,n)\Vert _{L^{1}(X,B)}
= \inf_{n
> 0} \sup_{ \Vert \lambda\Vert _{\infty,\mathcaligr{C}} \leq1} \frac{1}{n}
\int
_{X} \langle\lambda(x) , s_x(0,n) \rangle \,d\mu(x).
\]
We define the skew-product $\mathbb{Z}$-action $\hat{T}$ on the measurable
space $X \times\mathcaligr{C} $ with the product $\sigma$-algebra by
\[
\hat{T}_n(x,y) = (T_n x , c(n,x)^{*}.\lambda),\qquad x \in X, \lambda
\in\mathcaligr{C}.
\]
Note that if $n \geq0$ and $\lambda\in\mathcaligr{C}$, then
\[
\langle\lambda, s_x(0,n) \rangle= \sum_{k=0}^{n-1} \langle\lambda,
c(k,x).s_x(0,1) \rangle= \sum_{k=0}^{n-1} F(\hat{T}_k(x,\lambda)),
\]
where $F(x,\lambda) = \langle\lambda, s_x(0,1) \rangle$, and thus
\begin{eqnarray*}
A &=& \inf_{n > 0} \sup_{\Vert \lambda\Vert _{\infty,\mathcaligr{C}} \leq1}
\frac
{1}{n}\int_{X} \sum_{k=0}^{n-1} F(\hat{T}_k(x,\xi))\, d \delta
_{\lambda(x)}(\xi)\, d\mu(x) \\
&=& \inf_{n > 0} \sup_{\hat{\mu} \in M^1_\mu(X \times\mathcaligr{C})}
\frac{1}{n} \int_{X \times\mathcaligr{C}} \sum_{k=0}^{n-1} F(\hat
{T}_k(x,\xi))\,d\hat{\mu}(x,\xi),
\end{eqnarray*}
where $M_{\mu}(X \times\mathcaligr{C})$ denotes the space of probability
measures on $X \times\mathcaligr{C}$ which projects onto $\mu$ under
the canonical
map from $X \times\mathcaligr{C}$ to $X$. By standard disintegration
theory (see, e.g., \cite{Ar}), this space can be given a compact
metrizable topology arising
from the duality of $L^{1}(X,C(\mathcaligr{C}))$. Following the outline of
the proof in \cite{KaLe}, we take a sequence of elements $\hat{\mu
}_n$ in
$M^1(X\times\mathcaligr{C})$ such that
\[
\frac{1}{n} \int_{X} \sum_{k=0}^{n-1} F(T_k(x,\xi))\, d\hat{\mu
}_n(x,\xi) \geq A\qquad \forall n \geq1.
\]
This is possible due to the compactness of $M^1(X \times\mathcaligr
{C})$. Define
\[
\hat{\nu}_n = \frac{1}{n} \sum_{k=0}^{n-1} \hat{T}^k_*\hat{\mu}_n,\qquad
n \geq1.
\]
By sequential compactness, there exist a convergent subsequent and a
$\hat{T}$-invariant limit probability measure $\hat{\nu}_0$ in
$M^1_\mu
(X \times\mathcaligr{C})$ such that
\[
\int_{X \times\mathcaligr{C}} F(x,\xi) \,d\hat{\nu}_{0}(x,\xi) \geq A
\]
and thus the set of $\hat{T}$-invariant probability measures on $X
\times\mathcaligr{C}$ which project onto~$\mu$ and satisfy the above
inequality is a compact and convex
subset of $M^1_\mu(X \times\mathcaligr{C})$. By the Krein--Milman
theorem, there must be an extremal point $\nu$ in this set, and by a
standard argument (see, e.g., \cite{Ar}), this
point is an ergodic measure for $\hat{T}$. By Birkhoff's theorem and
the obvious inequality
\[
| \langle\xi, s_x(0,n) \rangle| \leq\Vert s_x(0,n)\Vert _{B}
\]
for all $\xi$ in the unit ball of $B^*$, we can conclude that
\[
A = \lim_{n \rightarrow\infty} \frac{1}{n} \Vert s_x(0,n)\Vert _{B} = \lim
_{n \rightarrow
\infty} \frac{1}{n} \langle\xi, s_x(0,n) \rangle
\]
for a co-null subset of $X \times\mathcaligr{C}$ with respect to the
measure $\nu$. If we assume that $(X,\mathfrak{F},\mu)$ is standard
Borel space, then we can use the Von Neumann
selection theorem (in complete analogy with \cite{KaLe}) and establish
the existence of a \textit{measurable} map $\xi\dvtx X \rightarrow
\mathcaligr{C}$
such that
\[
\lim_{n \rightarrow\infty} \frac{\langle\xi(x) , s_x(0,n) \rangle
}{n} = A
\]
for all $x$ in a co-null subset of $X$. We have established the
following theorem.
\begin{thm} \label{mult}
Suppose that $(X,\mathfrak{F},\mu)$ is a standard measure space with an
ergodic $\mathbb{Z}$-action. Suppose that $B$ is a separable Banach
space with
the Radon--Nikodym property and that $s$ is an
integrable Pettis cocycle with respect to a Borel cocycle $c$. Suppose
that there is a weak*-compact subset of $B_1^*$ such that
\[
\Vert s(0,n)\Vert _{1,B} = \sup_{\Vert \lambda\Vert _{\infty,B^{*}} \leq1} \int
_{X} \langle
\lambda(x) , s_x(0,n) \rangle\, d\mu(x)\qquad \forall n \geq1.
\]
There is then a measurable map $\xi\dvtx X \rightarrow\mathcaligr{C}$ such that
\[
\lim_{n \rightarrow\infty} \frac{1}{n} \langle\xi(x) , s_x(0,n)
\rangle= \lim_{n \rightarrow
\infty} \frac{1}{n} \int_{X} \Vert s_x(0,n)\Vert _{B}\, d\mu(x)
\]
almost everywhere on $(X,\mu)$.
\end{thm}

\begin{rema*}
The main reason for including the proof above is an application to
Kingman decompositions of subadditive cocycles which will be described
below. Note that the
restrictions on the Banach space $B$ and the measurability of $s$ are
fairly severe and exclude many interesting applications. For instance,
note that the case of
Pettis cocycles for $B = C(H)$, where $H$ is a compact metrizable
space, would generalize the celebrated multiplicative ergodic theorem
of Oseledec \cite{Os}.
In this situation, Theorem \ref{mult} was established for a certain
class of cocycles by Karlsson and Ledrappier in \cite{KaLe}. One
important feature of these cocycles
is an obvious choice of a sequence of weakly \textit{measurable} maps
$\eta_n \dvtx X \rightarrow B^{*}$ such that
\[
\Vert s_x(0,n)\Vert _B = \langle\eta_n(x) , s_x(0,n) \rangle
\]
for $n \in\mathbb{Z}$. This is no longer true for general cocycles in Banach
spaces. The Radon--Nikodym assumption on $B$ is a convenient way to
circumvent this problem.
\end{rema*}

An extension of Theorem \ref{mult} to conservative and ergodic actions
of $\mathbb{Z}$ on $\sigma$-finite measure spaces can be proven using
the same
techniques as in \cite{Bj},
where the Karlsson--Ledrappier ergodic theorem is extended to the
$\sigma
$-finite situation.

We now turn to the proof of an alternative Kingman decomposition for
random semimetrics induced by Pettis cocycles on reflexive and
separable Banach spaces. Let $\eta$
denote the disintegration of $\nu$ with respect to the canonical
projection onto measure $\mu$. Let $g \dvtx X \rightarrow\operatorname
{Isom}(B)$ be the
generator of the Borel cocycle $c$, that is,
\[
c(n,x) = g(x) \cdots g(T^{n-1}x)
\]
for $n \geq0$. For all $f \in L^1(X,B)$, we have
\begin{eqnarray*}
\langle\nu, \hat{T} f \rangle&=& \int_{X} \langle\eta(x) , g(x)
f(Tx) \rangle \,d\mu
(x) \\
&=& \int_{X} \langle(g(x))^{*} \eta(x) , f(Tx) \rangle\, d\mu(x) \\
&=& \int_{X} \langle(g(T^{-1}x)^{*}) \eta(T^{-1}x) , f(x) \rangle
\,d\mu(x) \\
&=& \int_{X} \langle\eta(x) , f(x) \rangle\, d\mu(x) = \langle\nu,
f \rangle.
\end{eqnarray*}
Thus, if $B$ is a reflexive Banach space, we conclude that
\[
\eta(Tx) = (g(x)^*)^{-1} \eta(x)
\]
or, equivalently,
\[
\eta(T^{k}x) = c(k,x)^* \eta(x)\qquad \forall k \geq1.
\]
Thus, we can rewrite the Birkhoff sum above as
\[
\langle\eta(x) , s_x(0,n) \rangle= \sum_{k=0}^{n-1} \langle\eta
(x) , c(k,x). f(T^k
x) \rangle= \sum_{k=0}^{n-1} \varphi(T^{k}x),
\]
where $\varphi(x) = \langle\eta(x) , f(x) \rangle$ satisfies
\[
\int_{X} \varphi(x)\, d\mu(x) = A.
\]
Furthermore, we obviously have
\[
\Vert s_x(0,n)\Vert _B \geq\sum_{k=0}^{n-1} \varphi(T^kx),\qquad n \geq1.
\]
We have proven the following weak version of Kingman's decomposition of
subadditive cocycles.
\begin{thm}[(Kingman decomposition)] \label{king}
Suppose that $s$ is an integrable Pettis cocycle with values in a
separable and reflexive Banach space, defined on a~standard probability
measure space with an ergodic $\mathbb{Z}$-action.
The random semimetric defined by
\[
\rho_x(m,n) = \Vert s_x(m,n)\Vert _B,\qquad  n,m \in\mathbb{Z},
\]
then decomposes as
\[
\rho_x(0,n) = \sum_{k=0}^{n-1} \varphi(T^k x) + r_n(x),
\]
where $\varphi$ is integrable on $(X,\mu)$ such that $\int_{X}
\varphi(x)\,d\mu(x)$ equals the drift of $\rho$ and~$r_n$ is a nonnegative
subadditive cocycle with drift $0$.
\end{thm}

\begin{rema*}
Kingman \cite{Ki} established a more general decomposition theorem for
integrable subadditive cocycles. Note, however, that Theorem \ref{king}
provides more information
about the decomposition. The restrictions on the measurability of $s$
and the Banach space $B$ in the theorem above seem to be necessary for
the methods described. However,
it is natural to ask for a canonical class $\mathcaligr{M}$ of Gelfand
cocycles on ergodic $G$-spaces and with values in Banach spaces with
separable pre-duals, such that for
any $s$ in $\mathcaligr{M}$ with values in $B$, there is an
$G$-equivariant and weakly*-measurable map $\eta\dvtx X \rightarrow B^*$
such that
\[
\Vert s_x(e,g)\Vert _{B} = p(g) \langle\eta_x , s_x(e,g) \rangle+ r_x(e,g),
\]
where $p \dvtx G \rightarrow\mathbb{R}$ is a weight function and $r_x$ is
negligible with
respect to $s$ in a~certain sense. In the case where $G = \mathbb
{Z}^d$ and
the seminorm $L$ in Theorem \ref{main2} is nondegenerate, this would
have interesting implications for generalized first passage percolation.
Indeed, this would imply a multiparameter version of Oseledec's theorem
with possible applications to infinite geodesics in random metric spaces.
\end{rema*}

\subsection{Horofunctions in random media} \label{horo}

Suppose that $\mathcaligr{H}$ is a separable Hilbert space and that $s\dvtx X \times\mathbb{Z}
^d \times\mathbb{Z}^d \rightarrow\mathcaligr{H}$ is a Bochner
cocycle in $L^{d,1}(X,\mathcaligr{H})$. Recall
that
\[
\rho_x(m,n) = \Vert s_x(m,n)\Vert _{\mathcaligr{H}},\qquad n,m \in\mathbb{Z}^d,
\]
defines a random semimetric on $\mathbb{Z}^d$. Suppose that $m$ is in
$\mathbb{Z}^d$
and define the horofunction at the point $m$, with respect to the
random semimetric $\rho$,
by
\[
h_m(n) = \rho(m,n) - \rho(m,0),\qquad n \in\mathbb{Z}^d.
\]
We want to study the behavior of $h_m$ as $m$ leaves finite subsets of
$\mathbb{Z}^d$. We will see that the limit exists along the sequence $m^{j}$
if and only if there is an element $\eta$ in the unit ball of $\ell
^1(\mathbb{Z}
^d)$ such that
\[
\lim_{j \rightarrow\infty} \frac{m_k^{j}}{|m|} = \eta_k,\qquad
k=1,\ldots,d.
\]
It will follow from the proof that the limit point is unique, that is,
independent of the particular sequence which converges to $\eta$. We
will denote
the limit point by~$h_\eta$ and refer to it as the \textit{horofunction
located at} $\eta$. Before we give the proof, we establish the
following simple
lemma.
\begin{lem} \label{lem}
Suppose that $m^{j}$ is a sequence in $\mathbb{Z}^d$ such that there
exists an
element $\eta$ in the unit ball of $\ell^{1}(\mathbb{R}^d)$, such that
$ m^{j}_k / |m| \rightarrow\eta_k $ for $k=1,\ldots,d$, where $|
\cdot| $
denotes the $\ell^1$-metric. Suppose that $s$ is a Bochner cocycle in
$L^1(X,\mathcaligr{H})$, where $\mathcaligr{H}$ is a separable
Hilbert space and $(X,\mu)$ is
an ergodic $\mathbb{Z}^d$-space. Then
\[
\lim_{m \rightarrow\eta} \frac{\Vert s_x(0,m)\Vert _\mathcaligr{H}}{|m|} =
\Biggl\Vert  \sum_{k=1}^{d} \eta
_k L_k \Biggr\Vert _{\mathcaligr{H}}
\]
almost everywhere on $X$ with respect to $\mu$. Here, $L_k = L(e_k)$,
$k=1,\ldots,d$, and $L$ is the continuous linear map in Theorem
\ref{thm:linear}. Conversely, the limit
\[
\lim_{m \rightarrow\eta} \frac{\Vert s_x(0,m)\Vert _\mathcaligr{H}}{|m|}
\]
exists almost everywhere on $X$ if and only if $m_{k}/|m|$ converges to
$\eta$.
\end{lem}

\begin{pf}
By Theorem \ref{thm:linear},
\[
\lim_{|m| \rightarrow\infty} \frac{\Vert  s_x(0,m) - L(m)\Vert _\mathcaligr
{H}}{|m|} = 0
\]
almost everywhere on $X$. Thus,
\[
\lim_{m \rightarrow\eta} \frac{\Vert s_x(0,m)\Vert _\mathcaligr{H}}{|m|} =
\lim_{m \rightarrow\eta}
\frac{\Vert s_x(0,m)-L(m) + L(m)\Vert _\mathcaligr{H}}{|m|} = \Vert  \eta_k L_k
\Vert _\mathcaligr{H}
\]
since $L(m) = \sum_{k=1}^{d} m_k L_k $ for all $m \in\mathbb{Z}^d$.
\end{pf}
In general, if $(Y,d)$ is a semimetric space, we define the
horofunction at a point~$y$ in $Y$ by
\[
h_y(y') = d(y,y') - d(y,0),\qquad  y' \in Y.
\]
If $d$ is a metric, the map $ y \mapsto h_y $ is injective.
Furthermore, if $(Y,d)$ is a proper metric space, that is, closed and
bounded sets
are compact, then the closure of the set $\{ h_y \}_{y \in Y}$ in
$C(Y)$ is compact by the Arzela--Ascoli theorem. In our case, the semimetric
$\rho$ is, in general, not a metric, nor is the topology it induces
proper. However, the notion of a horofunction is still well defined. We will
study the asymptotic behavior of the horofunctions with respect to the
random semimetric $\rho$ defined above in terms of $\mathcaligr{H}$-valued
cocycles. It turns out that a nice description is possible in this situation.

\begin{thm}
Suppose that $(X,\mu)$ is an ergodic $\mathbb{Z}^d$-space and $s\dvtx X
\times\mathbb{Z}^d
\times\mathbb{Z}^d \rightarrow\mathcaligr{H}$ is a Bochner cocycle
in $L^{d,1}(X,\mathcaligr{H})$, where $\mathcaligr{H}
$ is a real
separable Hilbert space. Let
\[
\rho(m,n) = \Vert s(m,n)\Vert _\mathcaligr{H},\qquad n,m \in\mathbb{Z}^d,
\]
denote the associated random semimetric on $\mathbb{Z}^d$. If $\eta$
is an
element in $\ell^1(\mathbb{R}^d)$
such that $\xi= \sum_{k=1}^{d} \eta_k L_k $ is a nontrivial element in
$\mathcaligr{H}$, where $L$ is the continuous linear map in Theorem
\ref{thm:linear},
then
\[
h_\eta(n) = \frac{2 \langle s(0,n),\xi\rangle}{\Vert \xi\Vert _\mathcaligr
{H}},\qquad n \in\mathbb{Z}^d,
\]
almost everywhere on $X$ with respect to $\mu$.
\end{thm}

\begin{pf}
The proof is a straightforward modification of the standard method for
computing horofunctions on a Hilbert space. If we suppose that
$s_x(0,n)$ and $s_x(0,m)$
are both nontrivial elements of $\mathcaligr{H}$, then
\begin{eqnarray*}
\Vert s_x(n,m)\Vert _\mathcaligr{H}- \Vert s_x(m,0)\Vert _\mathcaligr{H}&=& \frac
{\Vert s_x(n,0)+s_x(0,m)\Vert _\mathcaligr{H}^2
- \Vert s_x(0,m)\Vert _\mathcaligr{H}^2}{\Vert s_x(n,0)\Vert _\mathcaligr
{H}+\Vert s_x(m,0)\Vert _\mathcaligr{H}} \\
&=& \frac{\Vert s_x(n,0)\Vert _\mathcaligr{H}^2 + 2 \langle s_x(n,0),s_x(0,m)
\rangle_\mathcaligr{H}}{|m|}
\\
&&{}\cdot \frac{|m|}{\Vert s_x(n,0)\Vert _\mathcaligr
{H}+\Vert s_x(m,0)\Vert _\mathcaligr{H}}.
\end{eqnarray*}
By Lemma \ref{lem},
\[
\lim_{m \rightarrow\eta} \Vert s_x(n,m)\Vert _\mathcaligr{H}-
\Vert s_x(m,0)\Vert _\mathcaligr{H}= 2 \langle s_x(0,n),
\hat{\xi} \rangle_\mathcaligr{H}
\]
almost everywhere on $X$, where $\hat{\xi} = \xi/ \Vert \xi
\Vert _\mathcaligr{H}$.
\end{pf}

\begin{rema*}
It is still an open problem to compute the horofunctions at infinity
for the classical first passage percolation metrics. This would provide more
refined knowledge of the asymptotic geometry of these semimetric
spaces. It is expected that these horofunctions can be arbitrarily wild;
indeed, by a celebrated result of Meester and H\"aggstr\"om \cite
{HaMe}, essentially any convex shape in $\mathbb{R}^d$ can be obtained
as an
asymptotic shape of a classical first passage percolation generated by
ergodic $\mathbb{Z}^d$-actions.
\end{rema*}

\subsection{Reproducing kernel Hilbert spaces}

In this subsection, we will describe natural examples of Bochner
cocycles with values in separable Hilbert spaces. Let $(\mathcaligr
{H},K,o)$ be a
pointed reproducing Hilbert space. This means that $\mathcaligr{H}$ is
a Hilbert
space of measurable functions on a measurable space $(Y,\mathcaligr
{G})$ with a
fixed base point~$o$ in $Y$ and $K \dvtx Y \times Y \rightarrow\mathbb
{C}$ is a positive
definite reproducing kernel, that is, for all finitely supported
sequences $(c_i,y_i)$ in $\mathbb{C}\times Y$, we have the inequality
\[
\sum_{i,j} c_i \overline{c_j} K(y_i,y_j) \geq0
\]
and for all $y$ in $Y$, we have
\[
\langle K(y, \cdot) , f \rangle_\mathcaligr{H}= f(y)\qquad \forall f \in
\mathcaligr{H},
\]
where $\langle\cdot, \cdot\rangle_\mathcaligr{H}$ denotes the
inner product on $\mathcaligr{H}$.

Suppose that a locally compact group $G$ acts measurably on
$(Y,\mathcaligr{G})$.
In many cases, the action of $G$ on $Y$
can be lifted to an isometric action of $G$ on $\mathcaligr{H}$ so
that the
measurable metric
\[
d(y,y') = \Vert  K(y,\cdot) - K(y',\cdot) \Vert _\mathcaligr{H}^2,\qquad y,y' \in Y,
\]
is invariant under the action. Let $(X,\mathcaligr{F},\mu,T)$ be a
$\mathbb{Z}$-action and
suppose that $\pi$ is a unitary representation of $G$ on $\mathcaligr
{H}$. Given a
measurable
map $g \dvtx X \rightarrow G$, we define the isometry on
$L^2(X,\mathcaligr{H})$ by
\[
\hat{T}f(x) = \pi(g(x)).f(Tx),\qquad f \in L^2(X,\mathcaligr{H}),
\]
almost everywhere on $X$ and we let
\[
f^*(x) = K(g(x)o,\cdot) - K(o,\cdot),\qquad x \in X.
\]
Let $s$ be the Bochner cocycle generated by $f^*$ and the action $\hat
{T}$. We will describe a situation where $\pi$ can be chosen so that
\[
d(Z_n(x)o,o) = \Vert  s_x(0,n) \Vert ^2_\mathcaligr{H}\qquad \forall n \in\mathbb{Z},
\]
where $Z_n$ is the Borel cocycle generated by $g$ and $T$. We believe
that this is a fairly general phenomenon.

Let $\mathbb{D}$ be the Poincar\'e disc, that is, the unit disc in
$\mathbb{C}$ with
the distance function~$\beta$ given by
\[
\beta(o,z) = \log\frac{1+|z|}{1-|z|},\qquad z \in\mathbb{D},
\]
where $o$ is the origin, and extended to all pairs $(z,z')$ in $\mathbb{D}
\times\mathbb{D}$ by isometry. The isometry group $G$ of $\mathbb
{D}$ is isomorphic to
the M\"obius group $PSL_2(\mathbb{R})$. The large scale behavior of
$\beta$ can be equivalently described by the metric (see \cite{Wick}
for a more details)
\[
d(z,z') = \Vert  K(z,\cdot) - K(z',\cdot) \Vert _\mathcaligr{H}^2,
\]
where $(\mathcaligr{H},K)$ is the normalized Dirichlet reproducing
kernel Hilbert
space \cite{Wick} on~$\mathbb{D}$, that is, the reproducing Hilbert
space of
holomorphic functions $\phi$ on~$\mathbb{D}$
with $\phi(o) = 0$ and subject to the integrability condition
\[
\Vert \phi\Vert _\mathcaligr{H}= \biggl( \int_{\mathbb{D}} |\phi'(z)|^2\, dA(z)
\biggr)^{1/2} < \infty,
\]
where $A$ is the Euclidean area measure on $\mathbb{D}$ and
\[
K(z,z') = -\log( 1-z \overline{z'} ),\qquad (z,z') \in\mathbb{D}.
\]
The precise relation between the metrics $\beta$ and $d$ is discussed,
in a slightly different language, in the paper \cite{Wick}. In this
example, the
representation $\pi$ can be chosen to be
\[
\pi(g).\phi(z) = \phi(g^{-1}z) - \phi(o),\qquad z \in\mathbb{D}.
\]
A discussion about the relevance of the metric $\beta$ and the Borel
cocycle $Z$ to random Schr\"odinger equations can be found in \cite{KaLe}.

\subsection{Rates of convergence} \label{conv}

In this subsection, we will prove quantitative statements about the
convergence to a limit shape under certain conditions. Our results will
not apply to
classical first passage percolation, where deep results have been
established in a series of paper (see, e.g., \cite{Be,Kes,Zh}). We
will restrict the study to Bochner cocycles with values in Hilbert
spaces. This allows for certain spectral measure computations to be
performed and
the methods will not generalize beyond uniformly convex Banach spaces.
In particular, $L^{\infty}$-spaces, which would be the relevant spaces for
classical first passage percolation, are certainly out of reach.

Let $s$ denote a Bochner cocycle on a $\mathbb{Z}^d$-space $X$ with
values in a
Hilbert space~$\mathcaligr{H}$. By the additivity and equivariance
properties of
$s$, we
note that
\[
s_x(0,ne_1) = \sum_{k=0}^{n-1} s_x\bigl(ke_1,(k+1)e_1\bigr) = \sum_{k=0}^{n-1}
\lambda(k).s_{T_{ke_1}x}(0,e_1)\qquad \forall n \in\mathbb{Z}^d,
\]
where $\lambda$ is an isometric representation of $\mathbb{Z}^d$ on
$\mathcaligr{H}$. For
notational convenience, we define $f(x) = s_x(0,e_1)$. By standard
Hilbert space calculations,
we have
\begin{eqnarray*}
&&\frac{1}{n^2} \int_{X} \Vert s_x(0,ne_1)\Vert _\mathcaligr{H}^2 \,d\mu(x) \\
&&\qquad=
\frac
{1}{n^2} \sum_{j,k=0}^{n-1} \int_{X} \langle\lambda
(j).f(T_{ke_1}x), \lambda
(k).f(T_{ke_1}x) \rangle_\mathcaligr{H}\,d\mu(x) \\
&&\qquad= \frac{1}{n^2}\sum_{j,k=0}^{n-1} \bigl\langle f(x) , \lambda
(k-j).f\bigl(T_{(k-j)e_1}x\bigr) \bigr\rangle_\mathcaligr{H}\,d\mu(x) \\
&&\qquad= \sum_{k=-n}^{n} \frac{(n-|k|)}{n^2} \int_{X} \langle f(x) ,
\lambda
(k).f(T_{ke_1}x) \rangle_\mathcaligr{H}\,d\mu(x).
\end{eqnarray*}
We introduce the unitary operator $U$ on $L^2(X,\mathcaligr{H})$,
defined by $U^k
f(x) = \lambda(k).f(T^k x)$. We note that the calculations above
establish the following
proposition.
\begin{pro}
Let $U$ be the unitary operator defined above. Then
\[
\lim_{n \rightarrow\infty} \frac{1}{n^2} \sum_{k=-n}^{n} (n-|k|)
U^k f = Pf,
\]
where $P$ is the projection onto the space of $U$-invariant vectors in
$L^2(X,\mathcaligr{H})$.
\end{pro}

\begin{rema*}
It should be remarked that the proposition is true for any unitary
operator on $L^2(X,\mathcaligr{H})$. This is an immediate consequence
of Von
Neumann's mean ergodic theorem.
We included the calculation above for later reference.
\end{rema*}

A slight reformulation of the above proposition is contained in the
following lemma.
\begin{lem} \label{erg}
Suppose that $s$ is a Bochner cocycle on a $\mathbb{Z}^d$-space $X$
with values
in a Hilbert space $\mathcaligr{H}$. Then
\[
\lim_{n \rightarrow\infty}\frac{1}{n}
\Vert s(0,ne_1)\Vert _{L^2(X,\mathcaligr{H})} =
\Vert Ps(0,e_1)\Vert _{L^2(X,\mathcaligr{H})},
\]
where $P$ is the projection onto the space $U$-invariant vectors in
$L^2(X,\mathcaligr{H})$.
\end{lem}

Suppose that $\Vert f\Vert _{L^2(X,\mathcaligr{H})}=1$ and let $\nu_f$ denote the
probability measure on $\mathbb{T}$ such that $\hat{\nu}_f(n) =
\langle U^n f , f \rangle
$ for all $n \in\mathbb{Z}$.
We are interested in the asymptotic behavior of the sequence
\[
R_n = \Biggl\Vert \sum_{k=0}^{n-1} U^k f\Biggr\Vert _{L^2(X,\mathcaligr{H})}^2 - n^2
\Vert Pf(x)\Vert _{L^2(X,\mathcaligr{H})}^2.
\]
By Lemma \ref{erg}, we may assume that $Pf = 0$ in $L^2(X,\mathcaligr
{H})$. Several
papers have been written on the analogous situation in the case of
classical first passage percolation; see, for example, the papers \cite{Al93,Be} and \cite{Zh}. In our situation, we prove the following
analog of
Kesten's inequality in \cite{Kes}.
\begin{thm}
Suppose that $\nu_f$ is absolutely continuous with respect to the Haar
measure $m$ on $\mathbb{T}$ and that $ \frac{d\nu_f}{dm}$ is
continuous at $0$.
There then exists a constant~$C$ such that
\[
\Biggl| \Biggl\Vert \sum_{k=0}^{n-1} U^k f(x) \Biggr\Vert _{L^2(X,\mathcaligr{H})}^2 - n^2
\Vert Pf(x)\Vert _{L^2(X,\mathcaligr{H})}^2 \Biggr| \leq Cn\qquad\forall n \in\mathbb{N}.
\]
\end{thm}
\begin{pf}
By Lemma \ref{erg}, we can, without loss of generality, assume that \mbox{$Pf
= 0$} as an element of $L^2(X,\mathcaligr{H})$. Thus,
by the calculation above, we have
\begin{eqnarray*}
\frac{\Vert \sum_{k=0}^{n-1} U^k f\Vert_{L^2(X,H)}}{\sqrt{n}} &=& \biggl(\int
_{\mathbb{T}} \sum_{|k|\leq n} \biggl(1- \frac{|k|}{n} \biggr) e^{2\pi i k
\theta} \,d\nu_f(\theta)\biggr)^{1/2} \\
&=& \biggl(\int_{\mathbb{T}} F_n(\theta)\, d\nu_f(\theta)\biggr)^{1/2},
\end{eqnarray*}
where $F_n$ denotes the Fej\'er kernel. Thus, if $\frac{\nu_f}{dm}$ is
continuous at $0$, then the limit stays bounded for large $n$ and we
are done.
\end{pf}

\section*{Acknowledgments}
The author would like to thank Yves Derriennic and Anders Karlsson for
interesting discussions, and the anonymous referee for suggesting many
clarifying remarks.
%
\def\cprime{$'$}

\printaddresses

\end{document}